\documentclass[12pt,a4paper]{article}

\usepackage{amssymb}

\usepackage{amsmath}

\newtheorem{theorem}{Theorem}[section]
\newtheorem{lemma}{Lemma}[section]
\newtheorem{proposition}{Proposition}[section]

\newtheorem{corollary}{Corollary}[section]
\newtheorem{remark}{Remark}[section]

\usepackage{hyperref}
\usepackage{graphicx}
\graphicspath{ {./images/} }

\newcommand{\mr}{\mathbb R}
\newcommand{\mn}{\mathbb N}

\newcommand{\mc}{\mathbb C}

\newcommand{\SV}{{\cal V}}
\newcommand{\SK}{{\cal K}}
\newcommand{\SD}{{\cal D}}
\newcommand{\SU}{{\cal U}}
\newcommand{\SW}{{\cal W}}

\newcommand{\SM}{{\cal M}}

\newcommand{\be}{\begin{equation}}
\newcommand{\ee}{\end{equation}}
\newcommand{\RM}{{\rm M}}
\newcommand{\RF}{{\rm F}}

\newcommand{\gd}{\delta}

\newcommand{\pn}{\par  \noindent}
\newcommand{\bn}{\par \vspace {0.25cm} \pn}
\newcommand{\ot}{\otimes}

\newcommand{\SB}{{\cal B}}
\newcommand{\SA}{{\cal A}}

\newcommand{\GD}{{\Delta}}
\newcommand{\GO}{{\Omega}}
\newcommand{\ei}{{\bf 1}}
\newcommand{\rd}{\, {\rm d}}

\newcommand{\GG}{\Gamma}
\newcommand{\GL}{\Lambda}
\newcommand{\SE}{{\cal E}}
\newcommand{\SF}{{\cal F}}

\newcommand{\SH}{{\cal H}}
\newcommand{\RL}{{\rm L}}
\newcommand{\RP}{{\rm P}}
\newcommand{\RT}{{\rm T}}
\renewcommand{\ge}{\epsilon}

\newcommand{\id}{{\rm id}}
\newcommand{\SC}{{\cal C}}
\newcommand{\RU}{{\rm U}}
\newcommand{\gl}{\lambda}

\newcommand{\RE}{{\rm E}}
\newcommand{\ga}{\alpha}
\renewcommand{\gg}{\gamma}

\newcommand{\ri}{{\rm i}\,}

\newcommand{\kernel}{{\rm kern }\, }

\renewcommand{\phi}{\varphi}

\newcommand{\ol}{\overline}

\renewcommand{\mit}{{\, \vert \, }}

\newcommand{\norm}{{\vert \vert}}

\title{{\bf Non-commutative stochastic processes  with independent increments}} 
\author{Michael Sch\"urmann \\
Institut für Mathematik und Informatik \\
Universität Greifswald} 
\date{\today}

\begin{document}

\maketitle

\pn
This article is on the research of Wilhelm von Waldenfels in the mathematical field of quantum (or non-commutative) probability theory. Wilhelm von Waldenfels certainly was one of the pioneers of this field. His idea was to work with moments and to replace polynomials in commuting variables by free algebras which  play the role of algebras of polynomials in non-commuting quantities. Before he contributed to quantum probability he already worked with free algebras and free Lie algebras. One can imagine that this helped to create his own special algebraic method which proved to be so very fruitful. He came from physics. His PhD thesis, supervised by Heinz König, was in probability theory, in the more modern and more algebraic branch of probability theory on groups. Maybe the three, physics, abstract algebra and probability, must have been the best prerequisites to become a pioneer, even one of the  founders, of quantum probability. 
\par
We concentrate on a small part of the scientific work of Wilhelm von Waldenfels. The aspects of physics are practically not mentioned at all. There is nothing on his results  in classical probability on groups (Waldenfels operators). This is an attempt to show how the concepts of non-commutative notions of independence and of L\'evy processes on structures like Hopf algebras developed  from the ideas of Wilhelm von Waldenfels.

\section{Algebraic central limit theorem}
If nothing else is said, algebras are understood to be associative and over the field  of complex numbers.
\bn
Gaussian (or normal) distributions arise as \emph{central limits}.
A gaussian distribution is given by its covariance matrix $Q = \bigl ( Q_{ij} \bigr) _{i, j \in [d] }$ with $d \in \mn$, $Q_{ij} \in \mr$, and where  $[d]$ denotes the set $\{ 1, \ldots , d \}$.
Here $d$ is the degree of freedom of the distribution $\gg _Q$, and $Q_{ij} = \int _{\mr ^d} x_i x_j \rd \gg _ Q$. The (real) $d \times d$-matrix $Q$ is positive semi-definite. In particular, $Q$ is symmetric and we  have $Q_{ij} = Q_{ji}$. By definition, a $d$-dimensional real random variable $X : \GO \to \mr ^d$ over a probability space $(\GO , \SF , \RP )$ is $\gg _Q$-distributed if $\RP \bigl( X \in B \bigr) = \gg _Q (B)$ for all Borel subsets $B$ of $\mr ^d $ or equivalently, if
the expectation values $\RE \bigl( f \circ X \bigr) $ equal $\int _{\mr ^d } f \rd \gg _Q =: \gg _Q (f)$ for all bounded continuous functions $f$ on $\mr ^d $.
The \emph{moments} $\gg _Q (x_1 ^{l_1} \ldots x_d ^{l _d} )$, $l_1 , \ldots , l _d \in \mn _0 = \mn \cup \{ 0 \}$ of the gaussian $\gg _Q$ are given by
putting $x_1 ^{l_1} \ldots x_d ^{l _d} = x_ {i _1} \ldots x_{i_ n}$, $n = l_1 + \ldots + l_d$, and
\be
\label{momgau}
\gg _Q (x_ {i _1} \ldots x_{i_ n} ) = 
\left\{ \begin{array}{ll} 
0 & \mbox{if } n \mbox{ is odd}  \\
\sum Q_{S_1}  \ldots Q_{S_k} & \mbox{if } n = 2k \mbox{ is even}.
\end{array} \right.
\ee
where the sum is taken over all pair partitions $\{ S_1 , \ldots , S_k \}$ of 
$[ 2k ]$ (that is all partitions of $[ 2k]$ into sets with exactly two elements) and where $Q_ {S} = Q_{i _k , i_l }$ for $S = \{ k , l\}$, $k < l$. In the case of the normal distribution all moments exist and they determine the distribution. We can think of $\gg _Q$ as a linear functional on the (real, unital) algebra $\mr [ x_1 , \ldots , x_d ]$ of polynomials in the \emph{commuting} indeterminates $x_1, \ldots , x_d$. Then $\gg _Q (p) = \int _{\mr ^d } p \rd \gg _q$ for $p \in \mr [ x_1 , \ldots , x_d ]$. 
\par
The fundamental idea of the paper \cite{GivW} of N. Giri and W. von Waldenfels is to replace  $\mr [ x_1 , \ldots , x_d ]$ by the algebra  $\mc \langle x_1 , \ldots , x_d \rangle$ of polynomials with complex coefficients in \emph{non-commmuting} indeterminates $x_1 , \ldots , x_d$. Then for a (complex) $d \times d$-matrix $Q$ the linear functional $\gg _Q$ on  $\mc \langle x_1 , \ldots , x_d \rangle$, formally again defined by equation (\ref{momgau}), is called \emph{qaussian} with covariance matrix $Q$.
Notice that now  $Q$ is \emph{hermitian} and not automatically symmetric. In particular, we now have $\gg _Q (x_i x_j ) = Q_{ij} = \ol{Q} _{ji} = \ol{Q_{ij}} = \ol{\gg _Q (x_j x_i )}$ which means that it can happen that $\gg _Q$ does not vanish on the ideal in $\mc \langle x_1 , \ldots , x_d \rangle$ generated by the element $x_i x_j - x_j x_i$ and that $\gg _Q$ cannot be identified with a gaussian on $\mc [ x_1, \ldots , x_d ]$ or $\mr [ x_1, \ldots , x_d ]$.
For example, if $d = 2$ and 
$
Q = \frac{1}{2} \begin{pmatrix} 1 & \ri  \\ - \ri & 1  \end{pmatrix},
$
 then $Q$ is positive semi-definite. Here 
$\gg _Q (x_1 x_2 ) = \frac12 \ri \neq - \frac12 \ri = \gg _Q (x_2 x_1 ) $, and $\gg _Q (x_1 x_2 - x_2 x_1 ) = \ri$. It can be shown that $\gg _Q$ vanishes on the ideal in $\mc \langle x_1, \ldots , x_d \rangle$ generated by the elements $x_1 x_2 - x_2 x_1 - \ri$ (see Theorem 2 of \cite{GivW}; cf. Theorem \ref{ACLT} below). Thus $\gg _Q$ is a linear functional on the non-commutative *-algebra generated by two self-adjoint elements $q = x_1$ and $p = x_2$ with the canonical commutation relation between the canonical pair $(p, q)$ of momentum and position.
\bn
We state the classical Central Limit Theorem (CLT) in its elementary form.
let $X_1 , X_2 , \ldots $ be a sequence of, in the classical sense, independent identically distributed $d$-dimensional random variables such that the second moments exist. Moreover, let $X_n$ be centralized and put $Q_{ij} = \RE (X_i X_j )$, $i, j \in [d]$. Then $Q$ is real positive semi-definite. The CLT says that
\[
\frac{X_1 + \ldots + X_n }{\sqrt n}  \underset{n \to \infty}{\longrightarrow} \gg _Q  \  \mbox{ in law}
\]
which means that 
\be
\label{ClassCLT}
\int _{\mr ^d} f \bigl( \frac{X_1 + \ldots + X_n }{\sqrt n}\bigr)  \rd \RP \underset{n \to \infty}{\longrightarrow} \int _{\mr ^d} f \rd \gg _Q
\ee
for all bounded continuous functions on $\mr ^d$. Non-constant polynomials are not bounded, so a polynomial cannot serve as an $f$ in (\ref{ClassCLT}). The left hand side of (\ref{ClassCLT}) even makes no sense for a polynomial $f$, unless \emph{all} moments of $X_n$ exist.
If all moments exist, an \emph{algebraic} version of the CLT says that (\ref{ClassCLT}) holds for all polynomials $f \in \mc [ x_1, \ldots , x_d ]$. We may, for the moment, forget about positivity. Given a normalized, centralized (i.e. $\phi (\ei ) = 1$, $\phi (x_i ) = 0$) linear functional $\phi$ on $\mc [ x_1, \ldots , x_d ]$ the algebraic CLT 
\be
\label{AlgCLT}
\phi ^{\star n} \bigl( p( \frac{x_1}{\sqrt n} , \ldots , \frac{x_d }{\sqrt n} \bigr) ) \underset{n \to \infty}{\longrightarrow} \gg _Q (p)
\ee
holds  for all $ p \in \mc [ x_1, \ldots , x_d ]$. 
Here $ p( \frac{x_1}{\sqrt n} , \ldots , \frac{x_d }{\sqrt n})$ denotes the polynomial obtained from $p$ by replacement of $x_1 , \ldots , x_d$ by $\frac{x_1}{\sqrt n} , \ldots , \frac{x_d}{\sqrt n}$.   The statement (\ref{AlgCLT}) matches (\ref{ClassCLT}) for $f \in \mc [ x_1 , \ldots , x_d ]$, because  the distribution of the sum of $n$ independent identically distibuted random variables is the $n$-fold convolution of the distribution of $X_1$. Of course, we must define the convolution of (arbitrary) linear functionals on polynomial algebras which we will do now.
\par
Define the algebra homomorphism
\be
\label{conv}
\GD : \mc [ x_1 , \ldots , x_d ] \to  \mc [ x_1 , \ldots , x_d ]  \ot
\mc [ x_1 , \ldots , x_d ] 
\ee
by $\GD x_i = x_i \ot  \ei + \ei \ot x_i $ on the generators. The convolution product of two linear functionals $\phi$ and $\psi$ is then defined by $\phi \star \psi = (\phi \ot \psi ) \circ \GD$. 
\bn
In exactly the same manner the \emph{comultiplication} 
\be
\label{Conv}
\GD :  \mc \langle x_1 , \ldots , x_d \rangle \to
\mc \langle x_1 , \ldots , x_d \rangle  \ot
 \mc \langle x_1 , \ldots , x_d \rangle
\ee
can be defined in the non-commutative case. Again, for linear functionals $\phi$ and $\psi$, this time on $ \mc \langle x_1 , \ldots , x_d \rangle$, the convolution is given by $\phi \star \psi = (\phi \ot \psi ) \circ \GD$.  The algebraic CLT of  Wilhelm von Waldenfels reads as follows.
\begin{theorem}
\label{ACLT}
\text{\rm{(N. Giri, W. von Waldenfels 1978 \cite{GivW})}}
\pn
{\bf(a)}
Let $\phi$ be a normalized, centralized linear functional on $\mc \langle x_1 , \ldots , x_d \rangle$. Then (\ref{AlgCLT}) holds for all $ p \in \mc \langle x_1, \ldots , x_d \rangle$.
\pn
{\bf (b)} 
The linear functional $\gg _Q$ vanishes on the ideal in $\mc \langle x_1 , \ldots , x_d \rangle$ generated by the elements
\[
x_i x_j - x_j x_i - \bigl( Q_{ij} - Q_{ji}\bigr) , \  i, j \in [d] .
\]
\end{theorem}
{\it Proof}: In (\cite{GivW}) a proof based on direct combinatorial calculations is given. A proof using a convolution logarithm will be sketched in the next section. 
\bn
In fact, we will see that $\gg _Q$ is the convolution exponential $\exp _{\star} g _Q$ of the \emph{cumulant functional} $g _Q$ of $\gg _Q$ which is the linear functional on  $\mc \langle x_1, \ldots , x_d \rangle$ with
\be
\label{cumgau}
g _Q \bigl( x_{i _1} \ldots x_{i_n} \bigr) = 
\left\{ \begin{array}{ll}
0 & \mbox{if } n \neq 2 \\
Q_{i_1 , i_2} & \mbox{if } n = 2 
\end{array}  \right.
\ee  

\section{Schoenberg correspondence, tensor case}

Let us begin with a few remarks on tensor products of vector spaces. 
The tensor product of two vector spaces $\SV$ and $\SW$ in terms of bases is the vector space with basis $B \times C$ if $B$ and $C$ are bases of $\SV$ and $\SW$. This definition has the disadvantage that one has to use bases.
Alternatively and basis-free,  the tensor product is the vector space $\SV \ot \SW $ with basis $\SV \times \SW$ divided by the sub-vector space spanned by elements of the form $(\ga v_1 + v_2 ,  \beta w_1 + w_2 ) - \ga \beta (v_1 , w_1 ) - (v_2 , w_2 ) - \ga (v_1 , w_2 ) - \beta (v_2 , w_1 )$, $\ga , \beta \in \mc$, $v_1 , v_2 \in \SV$, $w_1 , w_2 \in \SW$.
For $v \in \SV$, $w \in \SW$ one writes $v \ot w$ for the equivalence class of $(v, w)$. Then every element of $\SV \ot \SW$ can be written in the form
\be
\label{tensor}
\sum_{i = 1} ^n v_i \ot w_i 
\ee
with suitable $n \in \mn$, $v_1 , \ldots , v_n \in \SV$ and $w_1 , \ldots , w_n \in \SW$. The representation (\ref{tensor}) of an element $\SV \ot \SW$ is not unique. A representation (\ref{tensor}) of an element in $\SV \ot \SW$ is called \emph{minimal} if $n$ is the minimum of all representations (\ref{tensor}) of this element.
 \begin{lemma}
 \text{{\rm (Paul-Andr\'e Meyer)}}
 \pn
 A representation (\ref{tensor}) of a non-zero element of $\SV \ot \SW$ is minimal if and only if $v_1 , \ldots , v_n$ are linearly independent in $\SV$ and $w_1 , \ldots , w_n $ are linearly independent in $\SW$.
 \end{lemma}
{\it Proof}: Each element $F \neq 0$ of $\SV \ot \SW$ has a representation 
(\ref{tensor}). If $n = 1$ we must have $v_1 \neq 0$ and $w_1 \neq 0$ because $F = v_1 \ot w_1 \neq 0$. -- Suppose that $F \neq 0$ has a minimal representation of form (\ref{tensor}) with $ n \geq 2$. If $v_1 , \ldots , v_n$ are not linearly independent, one of the $v_1 , \ldots  , v_n$ must be a linear combination of the others. By inserting this linear combination, a short computation shows that we obtain a representation (\ref{tensor}) of our element $F$ with $n - 1$ summands which is a contradiction to the minimality of $n$.  Similarly, $w_1 , \ldots , w_n$ have to be linearly independent.$\square$

\begin{lemma}
\label{inter}
Let $\bigl( \SV _i \bigr) _{i \in I}$ and $\bigl( \SW _i \bigr) _{i \in I}$
be two families of linear subspaces, $\SV _i $ of the vector space $\SV$ and  $\SW _i$ of the vector space $\SW$, indexed by the same index set $I$.
Then
\be
\bigcap _{i \in I} \SV _ i \ot \SW _ i = \bigl(\bigcap _{i \in I} \SV _ i \bigl) \ot   \bigl(\bigcap _{i \in I} \SW _ i \bigl)
\ee
\end{lemma}
{\it Proof}: The inclusion \lq$\supset$\rq \ is clear. --
For \lq$\subset$\rq \ we observe that $F \in \SV \ot \SW$ is in $\SV _i \ot \SW _i$ if and only if $(\phi \ot \id ) (F) \in \SW _i$ for all $\phi \in \SV '$ and $(\id \ot \psi ) (F) \in \SV _i$ for all $\psi \in \SW '$. 
Now for $F \in \bigcap _{i \in I} \SV _i \ot \SW _i$ we have $(\phi \ot \id ) (F) \in \SW _i$ and $(\id \ot \psi ) (F) \in \SV _i$ for all $i \in I$ which means $(\phi \ot \id ) (F) \in \bigcap _{i \in I} \SW _i$ and 
$(\id \ot \psi ) (F) \in \bigcap _{i \in I} \SV _i$.$\square$

\begin{lemma}
\label{name}
Let $\SU$, $\SV$ and $\SW$ be vector spaces and let 
\begin{equation*}
x \in \SU \ot \SV \ot \SW. 
\end{equation*}
Then there are $n, m  \in \mn$, linearly independent vectors $u_1 , \ldots , u_n \in \SU$, vectors $v_{i , j} \in \SV$, $i = 1, \ldots , n$, $j = 1, \ldots , m$, and linearly independent vectors $w_1 , \ldots , w_m \in \SW$ such that
\begin{equation}
\label{dreier}
x = \sum_{i = 1} ^n \sum_{j = 1} ^m u_i \ot v_{ij} \ot w_j .
\end{equation}
\end{lemma}
\pn
{\it Proof}:
It is easy to see that each element of $\SU \ot \SV \ot \SW$ can be written in the form 
\eqref{dreier} with $u_i \in \SU$, $v_{ij} \in \SV$ and $w_j \in \SW$. We must show that the $u_i$ and $w_j$ can be chosen in such a way that they are linearly independent.
Assume that in \eqref{dreier} $u_1$ is a linear combination of $u_2 , \ldots , u_n$ so that
$u_1 = \sum_{i = 2} ^n \ga _i u_i$ for some complex numbers $\ga _2 , \ldots , \ga _n$. 
Then
\begin{equation*}
x = \sum_{i = 2} ^n \sum_{j = 1}^m u_i \ot \bigl( \ga _i v_{1j} + v_{ij} \bigr) \ot w_j .
\end{equation*}
By eliminating the linearly dependent vectors $u_i$, we arrive at a representation \eqref{dreier} of $x$ with linearly independent $u_1 , \ldots , u_n$. 
In a similar manner the linearly dependent $w_j$ can be eliminated by changing the vectors $v_{ij}$.$\square$
\bn
A \emph{coalgebra} is a triplet $(\SC , \gd , \gd )$ consisting of a vector space $\SC$, a linear mapping $\GD : \SC \to \SC \ot \SC$ and a linear functional $\gd$ on $\SC$ such that the  condition $(\GD \ot \id ) \circ \GD = (\id \ot \GD ) \circ \GD$ of co-associativity and the counit property $(\gd \ot \id ) \circ \GD = \id = (\id \ot \gd ) \circ \GD$ are fulfilled. (We use the natural identification $\SC \ot \mc = \SC = \mc \ot \SC$.) In other words, a coalgebra is a co-monoid in the monoidal category of vector spaces with the above described tensor product of vector spaces. Sometimes we use the \emph{Sweedler notation} $\GD c = \sum c_{(1)} \ot c _{(2)}$.

\begin{corollary}
If $\SC _i$, $i \in I$, is a family of sub-coalgebras of the coalgebra $\SC$ then $\bigcap_{ \in I} \SC _i$ is again a sub-coalgebra of $\SC$.
\end{corollary}
\pn
{\it Proof}: By Lemma \ref{inter} $\GD \SC _i \subset \SC _i \ot \SC _i$ implies that $\GD \SC _i \subset (\bigcap_{i \in I} \SC _i ) \ot (\bigcap_{i \in I} \SC _i )$.$\square$
\bn
For a subset $\SM$ of a coalgebra $\SC$ we define the sub-coalgebra 
$\langle \SM \rangle$ generated by $\SM$ to be the intersection of all sub-coalgebras of $\SC$ which contain $\SM$.

\begin{theorem}
\label{FTC}
Let $\SM$ be a subset of a coalgebra $(\SC , \GD , \gd )$. Then
\begin{equation}
\# \SM < \infty \implies \dim \langle \SM \rangle < \infty
\end{equation}
\end{theorem}
{\it Proof}:
Since it is clear that the sum of two sub-coalgebras is again a sub-coalgebra, it suffices to proof that $\langle \{ c \} \rangle$ is finite-dimensional for each element $c \in \SC$.
\par
By Lemma \ref{name} we have 
\begin{equation}
\label{equation1}
(\GD \ot \id ) \GD c = \sum_{i = 1} ^n \sum_{j = 1} ^m u_i \ot v_{ij} \ot w_j
\end{equation}
with $u_1 , \ldots , u_n \in \SC$ linearly independent and $w_1 , \ldots , w_m \in \SC$ linearly independent.
Denote by $\SD$ the linear span  of the vectors $v_{ij}$, $i = 1, \ldots , n$, $j = 1, \ldots , m$ in $\SC$.
Then $\SD$ is a finite-dimensional linear subspace of $\SC$.
We will show that $\SD$ is a sub-coalgebra of $\SC$ containing $c$.
\par
In fact, we have
\begin{equation*}
c = (\gd \ot \id \ot \gd ) (\GD \ot \id ) \GD  
= \sum_{i = 1 } ^n \sum_{j = 1} ^m \gd (u_i ) \gd (w_j ) v_{ij} \in \SD .
\end{equation*}
Moreover,
\begin{align}
\label{coass2}
(\GD \ot \id \ot \id ) (\GD \ot \id ) \GD c
&= \sum_{j = 1} ^m \bigl( \sum_{i = 1} ^n (\GD u_i ) \ot v_{ij} \bigr) \ot w_j \\
&= (\id \ot \GD \ot \id ) (\GD \ot \id ) \GD c   \notag \\
&= \sum_{j = 1} ^m \bigl( \sum_{i = 1} ^n  u_i  \ot( \GD  v_{ij}) \bigr) \ot w_j \notag
\end{align}
Since $w_1 , \ldots , w_m$ are linearly independent, an application of $\id \ot \id \ot \id \ot \phi$ with $\phi \in \SC '$, $\phi (w_j ) = \gd _{jl}$ to \eqref{coass2} gives, for all $l = 1, \ldots , m$,
\[
\sum_{i = 1} ^n (\GD u_i ) \ot v_{il} = \sum_{i = 1} ^n u_i \ot (\GD v_{il})
\]
which is an element of $\SC \ot \SC \ot \SD$. 
Using $\phi \in \SC '$, $\phi (u_i ) = \gd _{ik}$, we have
\begin{align*}
& (\phi \ot \id \ot \id ) \bigl( \sum_{i = 1} ^n u_i \ot (\GD v_{il} ) \bigr)  \\
&\qquad = \GD v_{kl}  
= (\phi \ot \id \ot \id ) \bigl( \sum_{i = 1} ^n (\GD u_i ) \ot v_{il}  \bigr) \in \SC \ot \SD
\end{align*}
for all $k = 1, \ldots , n$, $l = 1, \ldots , m$.
The same argument, but now applying $\id \ot \id \ot \GD$ and $\id \ot \GD \ot \id$ to \eqref{equation1}, shows that $\GD v_{ij} \in \SD \ot \SC$.
Using again Lemma \ref{inter}, we finally obtain 
\begin{equation*}
\GD v_{ij} \in (\SC \ot \SD ) \cap (\SD \ot \SC ) = (\SC \cap \SD ) \ot (\SD \cap \SC ) = \SD \ot \SD .\square
\end{equation*}
Denote by $\RF (\SC )$ the system of finite dimensional sub-coalgebras of the coalgebra $\SC$. The system $\RF (\SC )$ is ordered by inclusion. If $\iota _{\SD , \SE }$, $\SD , \SE \in \RF (\SC )$, $\SD \subset \SE$ denote the inclusion maps, then $\bigl( \iota _{\SD, \SE } \bigr) _{\substack{\SD , \, \SE \in \RF (\SC ) \\ \SD \subset \SE }}$ is an inductive system.
Since for $c \in \SC$ we have $c \in \langle \{ c\} \rangle \in \RF (\SC )$, the coalgebra is the inductive limit of the system $\iota _{\SD , \SE }$. This fact can often be used to define a linear mapping on $\SC$ by defining it on finite-dimensional sub-coalgebras and by using compatibility with the $\iota _{\SD, \SE }$.
\par 
For instance, if $A : \SC \to \SC$ is linear and such that it leaves invariant all sub-coalgebras of $\SC$, we first define $\exp A$ on finite-dimensional sub-coalgebras of $\SC$, in the usual manner as an exponential on a finite-dimensional space. Next we observe that the restriction of  the exponential $\exp A $ on $\SE \supset \SD$ agrees with the exponential $\exp A$ on $\SD$. We define $\exp A$ on $\SC$ as the inductive limit. In particular, for $\psi \in \SC '$, the mapping 
$T _{\psi } = (\id \ot \psi ) \circ \GD$ leaves invariant all sub-coalgebras of $\SC$. We  define $\exp T_{\psi} $ on the whole of $\SC$ as an inductive limit by defining it on $\SD \in \RF (\SC )$ to be the exponential of the restriction of $T_{\psi} $ to $\SD$. Then 
\be
\bigl( \gd \circ \exp T_{\psi }\bigr) (c) = \lim_{n \to \infty} 
\sum_{k = 0} ^n \frac{\psi ^{\star k} (c)}{k !}
\ee
showing that the series $\sum_{k = 0} ^\infty \frac{\psi ^{\star k} (c)}{k !}$ converges for all $c \in \SC$. We denote the limit functional by $\exp _{\star} \psi$ and call it the \emph{convolution exponential} of $\psi$.
Again arguing with the inductive limit, we also have 
\be
\label{euler}
\bigl( \exp_{\star} \psi \bigr) (c) = \lim_{n \to \infty} \bigl( \gd + \frac{\psi}{n} \bigr) ^{\star n} (c) .
\ee
for all $c \in \SC$.

\begin{remark} {\rm We can now sketch a proof of the algebraic CLT \ref{ACLT} of von Waldenfels. The polynomial algebra $\mc \langle x_1 , \dots , x_d \rangle$ is a coalgebra with the comultiplication $\GD$ of (\ref{Conv}) and the counit $\gd$ defined as the algebra homomorphism $\gd : \mc \langle x_1 , \dots , x_d \rangle \to \mc$ with $\gd x_i = 0$, i.e. $\gd (p) = p(0)$. It is a nice calculation to see that $\gg _Q = \exp_{\star} g _ Q$. Next apply the \lq convolution logarithmic series\rq \ to the left of (\ref{AlgCLT}). This works again by an application of the inductive limit argument. Some estimation, just like for complex numbers (see \cite{Chung}, p. 184), shows that the $\log_{\star}$ of the left of (\ref{AlgCLT}) converges to $g _Q$ pointwise. From this it follows that $\phi ^{\star n} \bigl( p( \frac{x_1}{\sqrt n} , \ldots , \frac{x_d }{\sqrt n} \bigr) )$ converges to $\bigl( \exp _{\star} g_Q \bigr) (p) = \gg _Q (p)$.}
\end{remark}
Since $\gg _Q$ is a probability measure it is positive in the sense that $\gg _Q (p^* p ) \geq 0$ for all polynomials $p$. The involution $p \mapsto p ^*$ is given if we think of the indeterminates $x_i$ as \emph{self-adoint} elements that is $(x_i ) ^* = x_i$. However, in the non-commutative case we did not prove that $\gg _Q$ is positive. It is easy to check that $g_Q$ is \emph{conditionally positive}, which means  $g_Q (p ^* p ) \geq 0$ for all $p$ with $\gd ( p) = p (0) = 0$. \emph{Schoenberg correspondence} tells us that $\exp _{\star} \psi $ is positive if $\psi$ is conditionally positive and hermitian with $\psi (\ei ) = 0$. 
We are going to prove the Schoenberg correspondence for sesquilinear forms on coalgebras which is also fundamental  for the proof of the much more general Schoenberg correspondence of \cite{Malte}. It was  Wilhelm von Waldenfels who first proved Schoenberg correspondence for sesquilinear forms in the special case of anti-cocommutative coalgebras; see \cite{vW3}.
\par
The tensor product of two coalgebras $(\SC _i , \GD _i , \gd _i )$, $i = 1, 2$, is the coalgebra $(\SC _1  \ot \SC _2 , (\id \ot \tau \ot \id) \circ (\GD _1 \ot \GD _2 ), \gd _1 \ot \gd _2 )$ where $\tau : \SC_1 \ot \SC _2 \to \SC _2 \ot \SC _1$ denotes the flip operator. The complex conjugate coalgebra $(\ol \SC , \ol \GD , \ol \gd )$ of a coalgebra $(\SC , \GD , \gd )$ is the complex conjugate vector space of the vector space $\SC$, i.e. $\ol \SC = \{ \ol c \mit c \in \SC \}$ is a copy of $\SC$ with $\ol v + \ol w := \ol{v + w}$, $\gl \ol v = \ol{\ol \gl v }$. The comultiplication is given by $\ol \GD ( \ol c) = \sum \ol{c_{(1)}} \ot \ol{c_{(2)}}$ and the counit by $\ol \gd (\ol c ) = \ol{\gd (c)}$. A sesquilinear form $K$ on $\SC$ can be identified with a linear functional on $\ol \SC \ot \SC$ via $K(\ol c\ot d ) = K(c, d)$. 
\par
We define the convolution product $K \star L$ of two sesquilinear forms on a coalgebra to be the convolution of $K$ and $L$ with respect to the coalgebra structure of $\ol \SC \ot \SC$. A sesquilinear form $K$ is called positive semi-definite (or simply positive) if $ K(c,c)\geq 0$ for all $c \in \SC$, and it is called hermitian if $K (c , d ) = \ol{ K( d ,  c )}$ for all $c, d \in \SC$. A sesquilinear form $K$ on a coalgebra $\SC$ is called \emph{conditionally positive} if $K (c, c ) \geq 0$ for all $c \in \kernel \gd$ (that is for all $c \in \SC$ with $\gd (c) = 0$).
A graded version of the following is a generalization of the result \cite{vW3} of von Waldenfels to arbitrary graded coalgebras. We state and prove the result only in the non-graded setting, because in our applications to *-bialgebras a \lq symmetrization procedure\rq  \ always allows a reduction to the non-graded case.

\begin{theorem}
\text{{\rm (Schoenberg correspondence for sesquilinear forms)}}
\label{Schoenberg}
\pn
For a sesquilinear form $L$ on a coalgebra $\SC$ the following are equivalent:
\pn
{\bf (i)} $L$ is conditionally positive and hermitian. 
\pn
{\bf (ii)} $\exp_{\star} (t L ) $ is positive  for all $t \in \mr _+$.
\end{theorem}
{\it Proof}: We show that $\exp_{\star} L $ is positive if $L$ is conditionally positive and hermitian. In view of Theorem \ref{FTC} we may assume that $\SC$ is finite-dimensional. Let $L$ be conditionally positive, hermitian. We fix a a scalar product $\langle \ , \ \rangle : \SC \ot \SC \to \mc$ on $\SC$. Put $L_{\ge} := L + \ge \langle \ , \ \rangle$. We show that $\exp _{\star} L_{\ge}$ is positive for all $\ge > 0$. Then we finished, because 
\[
\bigl( \exp _{\star} L \bigr) (c, c) = \lim_{\ge \downarrow 0} \bigl( \exp _{\star} L_{\ge} \bigr) (c, c) \geq 0
\]
for all $c \in \SC$.
\par
By compactness of the unit ball $B = \{ c \in \SC \mit \norm c \norm \leq 1 \}$and because $L$ is conditionally positive, hermitian we have:
For all $\eta > 0$ there is a $\gg > 0$ such that
\be
\label{cond}
c \in B , \ \bigl\vert \gd (c) \bigr\vert ^2 \leq \gg \implies L (c, c ) > - \eta .
\ee
Fix $\ge > 0$. We claim that there is an $n _0 \in \mn$ such that $\bigl( \exp_{\star} \frac{L_{\ge}}{n_0}\bigr) (c, c ) \geq 0$ for all $c \in B$. Proof of this claim: let $\gg > 0$ such that (\ref{cond}) holds for $\eta = \ge$. For $c \in B$
\be
\label{estimate}
\bigl( \exp_{\star} \frac{L_{\ge}}{n}\bigr) (c, c ) = \bigl\vert \gd (c) \bigr\vert ^2 + \frac{L(c, c) + \ge }{n} + \frac{1}{n^2 } R_n(c)
\ee
where $\vert R_n (c) \vert \leq C$ for a suitable constant $C$. For $c \in B$ with 
 $\bigl\vert \gd (c) \bigr\vert ^2 \leq \gg$ we have  $L(c, c) + \ge > 0$. A short calculation shows that the convolution product of two hermitian sesquilinear forms is again hermitian. Thus $R_n (c)$ is real, and (\ref{estimate}) will be $\geq 0$ for $n$ large.
 \par
 Using the fact that the Schur product of two positive semi-definite matrices is again positive semi-definite, one shows that the convolution of two positive sesquilinear forms is again positive. We have proved that $ \exp_{\star} \frac{L_{\ge}}{n_0}$ is positive. Therefore, $\exp_{\star} L_{\ge} = \bigl( \exp_{\star} \frac{L_{\ge}}{n_0}\bigr) ^{\star n_0}$ is positive.$\square$
 \bn
We have seen that the polynomial algebras are examples of coalgebras. They are also equipped with an involution by assuming the indeterminates to be self-adjoint elements.
Then $\bigl( x_{i_1 } \ldots x_{i_n}\bigr ) ^* = x_{i_n } \ldots x_{i_1}$. Thus $\mc \langle x_1 , \dots , x_d \rangle$ (and also the commutative polynomial algebra) becomes a $*$-algebra. Since $\GD$ and $\gd$ are $*$-algebra homomorphisms (where the tensor product carries the usual involution), the polynomial algebras are 
Hopf $*$-algebras with the antipode given by extension of $x_i \mapsto - x_i$. The general definitions  read as follows: A coalgebra $(\SB , \GD , \gd )$ which is also an (associative, unital) algebra such that $\GD$ and $\gd$ are algebra homomorphisms is called a \emph{bialgebra}. An \emph{antipopde} $S$ on a bialgebra is a linear mapping $S : \SB \to \SB$ which satisfies the equations $M \circ(\id \ot S ) \circ \GD = \gd \, \ei = M \circ(S \ot \id  ) \circ \GD$. Here $M : \SB \ot \SB \to \SB$ is the multiplication map of $\SB$. A bialgebra which has an antipode is called a \emph{Hopf algebra}. A bialgebra which is also a $*$-algebra such that $\GD$ and $\gd$ are $*$-mappings is called a $*$-bialgebra. A $*$-bialgebra with antipode is called a Hopf $*$-algebra. 
A linear functional $\phi$ on  a $*$-algebra $\SA$ is called positive if $\phi (a ^* a ) \geq 0$ for all $a \in \SA$. If $\SA$ has a unit, a positive normalized linear functional is called a state. A linear functional $\psi$ on a $*$-bialgebra $\SB$ is called conditionally positive if $\psi (b ^* b ) \geq 0$ for all $b \in \kernel \gd$.
\par
A main application of Theorem \ref{Schoenberg} is

\begin{corollary}
\label{Schoenalg} 
For a linear functional $\psi$ on a $*$-bialgebra $\SB$ the following are equivalent:
\pn
{\bf (i)} $\psi$ is conditionally positive and hermitian with $\psi (\ei ) = 0$.
\pn
{\bf (ii)} $\exp_{\star} (t \psi )$ are states for all $t \in \mr _*$.
\end{corollary}
{\it Proof}: Apply Theorem \ref{Schoenberg} to $L (b, c ) := \psi (b ^* c )$.$\square$
\bn
In this article we are mainly interested in two classes of examples. The first are the  polynomial Hopf $*$-algebras which we will rewrite a bit: Let $\SV$ be a vector space which carries an anti-linear mapping $v \mapsto v ^*$ such that $\bigl( v ^* \bigr) ^* = v$. We form the tensor algebra
\[
\RT (\SV ) := \bigoplus_{n = 0} ^\infty \SV ^{\ot n}
\]
with the multiplication given by $(v_1 \ot \ldots \ot v_n ) (w_1 \ot \ldots \ot w_m ) := v_1 \ot \ldots v_n \ot w_1 \ot \ldots \ot w_m$. (If $\SV$ is $d$-dimensional and $x_1 , \ldots , x_d$ forms a basis of self-adjoint elements of $\SV$, we may identify $\RT (\SV )$ with $\mc \langle x_1 , \ldots , x_d \rangle$.) By extension of $*$ as an involution from $\SV \subset \RT (\SV )$ to  the algebra  $\RT (\SV )$  in the only possible way we define an involution on $\RT (\SV )$. Moreover, $\RT (\SV )$ becomes a Hopf $*$-algebra by putting $\GD v = v \ot \ei + \ei \ot v$, $\gd v = 0$ and $Sv = - v$.
\par
The second class of examples is inspired by the coefficient Hopf algebra $\SK [d]$ of the group $\RU _d$ of unitary $d \times d$-matrices. $\SK [d]$ is the sub-algebra of the algebra of continouus complex-valued functions on $\RU _d$ which is formed by the coefficients of the continuous finite-dimensional representations of $\RU _d$.
By a result of Hermann Weyl $\SK [d]$ is $*$-isomorphic to the algebra $\mc \bigl[ x _{ij} , x_{ij} ^* ; i, j  \in [d] \bigr] / J$ where $J$ denotes the ideal gemerated by the elements
\be
\label{unitaryRel}
\sum_{n = 1} ^d x_{in} x_{jn} ^*  -  \gd _{ij} \ei ,  \
\sum_{n = 1} ^d x_{ni} ^* x_{nj}   -  \gd _{ij} \ei , \ i, j \in [d] .
\ee
The $*$-isomorphism is given by $x_{ij }  \mapsto u_{ij}$ with $u_{ij} : \RU _d \to \mc$, $ u_{ij} (U) = U_{ij}$. By extension of $\GD _{ij} = \sum_{n = 1} ^d x_{in}\ot x_{nj}$, $\gd (x_{ij}) = \gd _{ij}$ and $S x_{ij} = x_{ji} ^*$ the $*$-algebra $\SK [d]$ becomes a Hopf $*$-algebra.
\par 
Wilhelm von Waldenfels introduced a \emph{non-commutative version} of $\SK [d ]$. He was studying a mathematical model for the light emission and absorption in a laser; see \cite{vW84}. In this model the time evolution  is described by unitary operators on  $\mc ^d \ot \SH$, $\SH$ an infinite-dimensional Hilbert space. Here the physical system $\mc ^d$ is coupled to a \lq heat bath\rq \ which is represented by $\SH$. The mappings $\xi _{ij} : \RU (\mc ^d \ot \SH ) \to \SB (\SH )$ from the group $\RU (\mc ^d \ot \SH )$  of unitary operators on $\mc ^d \ot \SH$ to the $C^*$-algebra $\SB (\SH ))$ of bounded operators on $\SH$ are defined by $\xi _{ij} (U) = U_{ij}$ for $U \in \RU (\mc ^d \ot \SH )$, $U$ considered as a $d \times d$-matrix with entries in $\SB (\SH )$. Denote by $\SK \langle d \rangle$ the sub-$*$-algebra generated by the $\xi _{ij} $ of the $*$-algebra of mappings from $\RU (\mc ^d \ot \SH )$ to $\SB (\SH )$ (which is a $*$-algebra with the pointwise structure coming from $\SB (\SH )$).
On the other hand, von Waldenfels considered the algebra $\mc \bigl\langle  x _{ij} , x_{ij} ^* ; i, j  \in [d] \bigr\rangle /\tilde J$ where $\tilde J$ is  the ideal, now in the non-commutative polynomial algebra, generated by formally the same elements (\ref{unitaryRel}). The following analogue of Hermann Weyl's classical result was found by von Waldenfels in the 1980's and published in \cite{GlvW} as a cooperation8 with Peter Glockner.

\begin{theorem}
There is a unique $*$-algebra homomorphism 
\[
\Phi : \mc \bigl\langle  x _{ij} , x_{ij} ^* ; i, j  \in [d] \bigr\rangle /\tilde J
\to \SK \langle d \rangle
\]
such that $\Phi \bigl( x_{ij} + \tilde J \bigr) = \xi _{ij}$.
\end{theorem}
See \cite{GlvW} for a remarkable proof which uses Hermann Weyl's  result for all $\RU _m$, $m \in \mn$. 
\par
In fact, it is not difficult to see that $\GD x_{ij} = \sum_{n = 0} ^d x_{in} \ot x_{nj}$ and $\gd x_{ij} = \gd _{ij}$ defines a $*$-bialgebra structure on $ \mc \bigl\langle  x _{ij} , x_{ij} ^* ; i, j  \in [d] \bigr\rangle$. Moreover, $\tilde J$ is a $*$-bi-ideal of this $*$-bialgebra so that $\SK \langle d \rangle$ itself is turned into a $*$-bialgebra. We sometimes write $x_{ij}$ for the equivalence class $x_{ij} + \tilde J$. 
\par
The $*$-algebra $\SK \langle d \rangle$ has been introduced by L.G. Brown in \cite{Brown}  and is now called the \emph{Brown-Glockner-Waldenfels algebra}. One shows that $\SK \langle d \rangle$ is not a Hopf $*$-algebra, i.e. it does not allow for an antipode!

\section{L\'evy processes, tensor case}

A good source for this section is \cite{Fra06} or \cite{Lachs1} or \cite{Ger1}.
\par
A stochastic process $X_t : \GO \to G$, $t\in \mr _+$, with values in a topological group $G$ is called a \emph{L\'evy process} on $G$ (over the probability space $(\GO , \SF , \RP )$) if the increment process $\bigl( X_{st} \bigr) _{0 \leq s \leq t}$, $X_{st} (g) := X_s (g ) ^{-1} X_t (g)$ fulfils the following conditions.
The distribution of the increment $X_{st}$ only depends on the difference $t - s$ (that is the distribution of $X_{st}$ equals the distribution of $X_{0, t - s}$), the increments $X_{t_1 , t_2} , \ldots , X_{t_n , t_{n + 1}}$ for consecutive times $0 \leq t_1 \leq t_2 \leq \ldots , t_n \leq t_{n + 1 }$ are stochastically independent, and $X_{st}$ converges weekly to the constant random variable $e = X_{s s}$ for $t \downarrow s$. 
\par
In order to avoid inverses one starts right away from an \lq increment process\rq \ $\bigl( X_{st} \bigr) _{0 \leq s \leq t}$ with the extra increment conditions $X_{rs } X_{st} = X_{rt}$, $X_{tt} = e$. Thus we may define a
 L\'evy process on a topological monoid $G$ to be a stochastic process $X_{st} : \GO \to G$, $0 \leq s \leq t$, which satisfies the increment condition and with stationary, independent increments in the above sense. In the non-commutative  case, where Hopf algebras will replace groups, this will mean that we will have L\'evy processes on $*$-bialgebras rather than Hopf $*$-algebras; see  the example of the $*$-bialgebra $\SK \langle d \rangle$ (Brown-Glockner-Waldenfels algebra) which is not a Hopf algebra. Since we wish to enter the non-commutative world we will have to \lq reverse all the arrows\rq \ and replace the measurable spaces by $*$-algebras of functions and later by not necessarily commutative $*$-algebras. 
\par 
The underlying probability space is, for example, replaced by the pair $\bigl( \RL ^\infty (\GO ), E )$ with $\RL ^\infty (\GO )$ the $*$-algebra of complex-valued bounded measurable function on $\GO$ and $\RE (F)$ the expectation value of $F \in \RL ^\infty (\GO )$. For example, if $G = \RU _d$ we could replace $\RU _d$ by its coefficient algebra $\SK [ d]$. The process becomes the family $j_{st} : \SK [d] \to \RL ^\infty (\GO )$ of $*$-algebra homomorphisms defined by $j_{st} (f) = f \circ X_{st}$, $f \in \SK [d]$. The increment condition now reads $j_{rs } \star j_{st}\  ( := M \circ (j_{rs} \ot j_{st} ) \circ \GD) = j_{rt}$ and $j_{tt} = \gd \, \ei$.
 Stationarity becomes $\RE \circ j_{st} = \RE \circ j_{0, t - s}$, week continuity now reads $\RE \circ j_{st}  \underset{t \downarrow s}{\longrightarrow} \gd$ pointwise, and independence of increments means that the state $\RE$ factorizes on the image algebras $j_{st} \bigl( \SK [ d ] \bigr)$, i.e.
\[
\RE \bigl( j_{t_1 , t_2 }(f_1 )  \ldots j_{t_n , t_{n + 1} } (f_n ) \bigr) 
=  \RE \bigl( j_{t_1 , t_2 }(f_1 ) \bigr) \ldots \RE \bigl( j_{t_n , t_{n + 1} }(f_n )\bigr) 
\]
for all $t_1 \leq  \ldots \leq t_{n + 1}$, $f_1 , \ldots , f_n \in \SK [ d]$.  In his work on light emission and absorption \cite{vW84} Wilhelm von Waldenfels replaced $\SK [d]$ by $\SK  \langle d \rangle$; see §2 of \cite{vW84} where von Waldenfels uses the notions $K \bigl( U (d)\bigr)$ for $\SK [d]$ and $(\mbox{handwritten } K ) \bigl(  (U (d) \bigr)$ for $\SK \langle d \rangle$. We will return to this work of von Waldenfels later.
\bn
What is a L\'evy process on a $*$-bialgebra over a \lq quantum probability space\rq ? The latter is a pair $(\SA , \Phi )$ where $\SA$ is a (unital) $*$-algebra and $\Phi : \SA \to \mc$ is a state on $\SA$. In the classical case one may choose $\SA = \RL ^\infty (\GO )$ and $\Phi = \RE$. A tensor L\'evy process on a $*$-bialgebra $\SB$ over $(\SA , \Phi )$ is a family $j_{st} : \SB \to \SA$, $0 \leq s \leq t$, of $*$-algebra homomorphisms such that $j_{rs} \star j_{st} = j_{rt}$ and $j_{tt} = \gd \, \ei _{\SA}$, and such that $\Phi \circ j_{st} = \Phi \circ j_{0 , t - s}$ and $\Phi \circ j_{st} \underset{t \downarrow s}{\longrightarrow} \gd$ and the sub-algebras $j_{t _1 , t_2 } (\SB ) , \ldots , j_{t_n , t_{n + 1}} (\SB )$ are \lq 
independent\rq  \ for all $t_1 \leq \ldots \leq t_{n + 1}$. 
\par
Here we use the notion of what nowadays sometimes is called \emph{tensor independence} to replace the classical stochastic independence. Tensor independence or, more generally, \emph{graded} tensor independence, goes back to the early work of another pioneer of quantum probability, Robin L. Hudson, see \cite{Hud} and \cite{CuHu}. Tensor independence would probably be the first answer of a physicist if he was asked for a reasonable notion of quantum or non-commutative independence: The sub-algebras $\SA _1 , \ldots , \SA _n$ are called \emph{tensor independent} in the state $\Phi$ if $\SA _k $ and $\SA _l$ commute for different values of $k, l \in [n]$ (that is $a_k a_l - a_l a_k = 0$ for all $a_k \in \SA _k$ and $a_l \in \SA _l$) \emph{and} if $\Phi (a_1 \ldots a_n ) = \Phi (a_1 ) \ldots \Phi (a_n )$ for all $a_1 \in \SA _1, \ldots , a_n \in \SA _n$.

\begin{proposition}
Let $j_{st}$ be a tensor L\'evy process on a $*$-bialgebra $\SB$. Then we have with $\phi _t := \Phi \circ j_{0, t}$, $t \in \mr _+$,
\begin{align}
\label{c1}
\phi _{s + t} &= \phi _s \star \phi _t \mbox{ for } s, t \in \mr _+  \\
\label{c2}
\phi _t (b) &\underset{t \downarrow s}{\longrightarrow} \gd (b)  \mbox{ for all } b \in \SB .
\end{align}
Moreover,  $\frac{1}{t} \bigl( \phi _t - \gd \bigr) (b)$ converges for $t \downarrow 0$  for all $b \in \SB$ to a limit  which we denote by $\psi (b)$. The linear functional $\psi$ on $\SB$ is conditionally positive and hermitian with $\psi (\ei ) = 0$. We have $\phi _t = \exp _{\star} \bigl( t \psi \bigr)$.
\end{proposition}
{\it Proof}:
It follows straightforward from the properties of a tensor L\'evy process that $\phi _t$ is a weekly continuous convolution semi-groups of states on $\SB$, in particular the validity of (\ref{c1}) and (\ref{c2}).
\par
For a $d$-dimensional, $d \in \mn$, sub-coalgebra $\SD$ of $\SB$ the linear operators $A_t$ with $A_t (b) = \bigl( \id \ot \phi _t \bigr) \circ \GD (b)$, $b \in \SD$, satisfy $A_{s + t} = A_s \circ A_t$ and $A_t \underset{t \downarrow 0}{\longrightarrow} \id _ {\SD}$. By a well-known result for matrices $A_t$ must be of the form $A_t = \exp \bigl( t G\bigr)$ with $G$ a $d \times d$-matrix and $\frac{1}{t} \bigl( A_t - \id _{\SD} \bigr) \underset{t \downarrow 0}{\longrightarrow} G$. Now everything follows from the fact that the coalgebra $\SB$ is the inductive limit of its finite-dimensional sub-coalgebras; see Section 2.$\square$
\bn
The conditionally positive linear functional of the proposition is called the \emph{generator} of the tensor L\'evy process $j_{st}$. The generator determines the numbers $\Phi \bigl( j_{s _1 , t_1} (b_1 )  \ldots j_{s_n , t_n } (b_n ) \bigr)$ where $n$ runs through the natural numbers, $s_i , t_i \in \mr _+$, $b_1 , \ldots , b_n \in \SB$, which are the numbers relevant for the process. In this sense, $\psi$ determines the process. The converse construction is now important. Start from a conditionally positive, hermitian linear functional $\psi$ on $\SB$. Then by Schoenberg correspondence (Corollary \ref{Schoenalg}) $\phi _t = \exp_{\star} \bigl( t \psi \bigr)$ is a convolution semi-group of states on $\SB$. Inspired by the Daniel-Kolmogorov construction for stochastic processes of classical probability theory, in the 1ate 1980's Luigi Accardi achieved, in co-operation with Wilhelm von Waldenfels and M. Schürmann, the construction of a tensor L\'evy process with a given convolution semi-group as an inductive limit; see \cite{waldenfels88}. This established a 1-1-correspondence between conditionally positive, hermitian linear functionals on a $*$-bialgebra $\SB$ (the generators) and (equivalence classes) of tensor L\'evy processes on the $*$-bialgebra $\SB$.
\bn
Let us turn to a more detailed investigation of tensor L\'evy processes on  Hopf $*$-algebra of type $\RT (\SV )$ with $\SV$ a vector  space with involution $*$. Generators on $\RT (\SV )$ can be described as follows. Let $H$ be a pre-Hilbert space, i.e. a complex vector space with a fixed scalar product $\langle \ , \ \rangle$. Denote by $L _a (H)$ the $*$-algebra of adjointable linear operators on $H$. (A linear $T : H \to H$ is by definition in $L_a (H)$ if there is a linear $T^* : H \to H$ such that $\langle x , T y \rangle = \langle T ^* x , y \rangle$ for all $x, y \in H$.)
Let $\bigl( \rho _0 , \eta _0 , \psi _0 \bigr)$ be a triplet formed  by a linear $*$-map $\rho _0 : \SV \to L_a (H)$, a linear map $\eta _0 : \SV \to H$ and a hermitian linear functional $\psi _0$ on $\SV$. Then an interesting calculation shows that there is a unique generator $\psi$ on $\RT (\SV )$ such that $\psi _0$ is the restriction of $\psi$ to $\SV \subset \RT (\SV )$, for $v_1 , v_2 \in \SV$ we have $\psi (v_1 \ot v_2 ) = \langle \eta _0 (v_1 ^* ) , \eta _0 (v_2 ) \rangle$, and
\[
\psi (v_1 \ot \ldots \ot v_n ) = 
\langle \eta _0 (v_1 ^* ) , \rho _0 (v_2 )  \ldots \rho _0 (v_{n - 1}) \eta _0 (v_n ) \rangle
\]
for $ n \geq 3$ and $v_1 , \ldots, v_n \in \SV$. This construction yields a 1-1-correspondence between generators and triplets of this kind.
\par
We form the pre-Hilbert space $\GG (H) := \bigoplus _{n = 0} ^\infty K ^{\ot _s n}$ with $K = \RL ^2 (\mr _+ ) \ot H$ where $K ^{\ot _s n}$ denotes the symmetric tensor product of $n$ copies of $K$. On this \emph{symmetric Fock space} $\GG (H)$ we define the operators 
\begin{align} 
A_t (x) &= A (\chi _{[0, t[} \ot x)  \\
A_t ^* (x) &= A (\chi _{[0, t[} \ot x)^*  \\
\GL _ t (T) &= \GL (\chi _{[0, t]} \ot T )
\end{align}
for $x \in H$, $T \in L_a (H)$. Here $A (\xi )$, $A (\xi ) ^*$, $\xi \in K$ and $\GL ( \RT )$, $T\in L_a (K)$ denote the annihilation, creation and additive second quantization (= preservation) operators on symmetric Fock space respectively; cf. \cite{Par92}. Let $\Phi$ be the state on $L_a \bigl( \GG (H)\bigr)$ given by the unit vector $(1, 0 , 0, \ldots )$ (called \lq vacuum\rq ). Consider the operators ($v \in \SV$, $ t \in \mr _+$) 
\be
\label{additive}
F _t (v ) = A _t \bigl( \eta _0 (v ^* )\bigr) + \GL _t \bigl( \rho _0 (v )\bigr) + A _t ^* \bigl( \eta _0 (v )\bigr) + \psi _0 (v) \, t
\ee
in $L_a (\GG (H) \bigr)$. Define homomorphisms $j_{st} : \RT (\SV ) \to L_a (\bigl( \GG (H)\bigr)$ as the extensions of $v \mapsto F _t (v) -  F_s (v)$ to $\RT (\SV )$.
Then some standard argumentations show that 
$j_{st}$ is a realization of a tensor L\'evy process with generator $\psi$ given by the triplet $(\rho _0 , \eta _0 , \psi _0 )$. The underlying quantum probability space is $\bigl (L _a \bigl( \GG (H)\bigr), \Phi \bigr)$. We have
\begin{theorem}
Let $\psi$ be a conditionally positive, hermitian linear functional on $\RT (\SV )$  with pre-Hilbert space $H$ and triplet $(\rho _0 , \eta _ 0 , \psi _0 )$ as above. Then equation (\ref{additive}) defines a tensor L\'evy process on $\RT (\SV )$ with generator $\psi$.
\end{theorem}
As an example, consider the linear functional $g _Q$ defined by equation(\ref{cumgau}). We assume that $Q$ is positive semi-definite. For simplicity let us also assume that $Q$ is regular, i.e. that $Q$ is positive definite. It is clear that $g _Q$ is hermitian and conditionally positive so that it is the generator of a tensor L\'evy process on $\mc \langle x_ 1 , \ldots , x_d \rangle$. 
An interesting calculation shows that indeed $\gg _Q = \exp_{\star} g _Q$!
The pre-Hilbert space $H$ in this case is $\mc ^d$ equipped with the scalar product $\langle x, Q y \rangle$. The map $\eta$ vanishes on all components $\SV ^{\ot n}$ with $n \neq 1$ where we put $\SV = \mc ^d = H$ equal to the linear span of $x_1 , \ldots , x_d$ and use the identification of $\mc \langle x_ 1 , \ldots , x_d \rangle$ with $\RT (\SV )$. Now for $v \in \SV$ we have $\eta( v)  = v \in H$. Moreover, $\rho (v ) = 0$ for all $v \in \SV$. The L\'evy process is given by
\[
F_t (x_ i ) = A_t (x_i ) + A_t ^* (x_i ),  \    i = 1, \ldots , d.
\]
We obtain $d$ realizations of Brownian motions as operator processes on the same symmetric Fock space. Of course, in general  the $d$  Brownian motions do not commute. We arrive at Robin Hudson's \emph{quantum Wiener process}.
\bn
Let us have a look at the general tensor case. A generator $\psi$ of a tensor L\'evy process $j_{st}$ on a $*$-bialgebra $\SB$ is associated with a pre-Hilbert space $H$, a $*$-algebra map $\rho : \SB \to L_a (H)$ and a linear surjective linear map $\eta : \SB \to H$ via the cohomological equations
($b, c \in \kernel \gd$)
\begin{align*}
\rho (b) \eta (c) &= \eta (bc)\\
 \langle \eta (b ^* ) , \eta (c ) \rangle &=   \psi (b c)  ;
\end{align*}
see \cite{MSchue91b} and \cite{Fra06}. Put $\SB _0 = \kernel \gd$. Then use the restrictions $\rho _0 := \rho \restriction \SB _0$, $\eta _0 = \eta \restriction \SB _0$ and $\psi _0 = \psi \restriction \SB _0$ as a triplet to define the generator of a tensor L\'evy process $k_{st}$ on $\RT (\SB _0 )$ which we call the \emph{generator} process of the tensor L\'evy process $j_{st}$ on $\SB$. It will be given by $F_t $ of (\ref{additive}). Since $\RT (\SB _0 )$ is a Hopf $*$-algebra, the process $k_{st}$ is determined by the process $k_t$ with $k_t := k _{0t}$.
\par
Now we use \emph{quantum stochastic calculus} as it was developed by Robin L. Hudson and K.R. Pathasarathy in the 1980's; see \cite{HudPar} and \cite{Par92}.
The tensor L\'evy process $j_{st}$ with generator $\psi$ can always be realized on the symmetric Fock space $\GG (H)$ as the solution of the quantum stochastic differential equation (QDE) 
\be
\label{QSDE}
\rd j_{st} = j_{st} \star \rd k_t , \  j_{ss} =\gd
\ee
which is driven by the generator process $k_t$;
see \cite{MSchue91b} and \cite{Lachs1}) for details. (Here we identify $\SB$ with $\mc \ei \oplus \SB _0 \subset \RT (\SB _0 )$.)
 Again there is also pioneering work of Wilhelm von Waldenfels in this field of quantum probability. See again his paper \cite{vW84} where he constructed a tensor
  L\'evy process on $\SK \langle d \rangle$ using a multiplicative quantum stochastic It\^{o} integral.
\bn
The process which von Waldenfels constructed in \cite{vW84} is a tensor L\'evy process on the $*$-bialgebra $\SK \langle d \rangle$. 
\par
Given a tensor L\'evy process $j_{st}$ on $\SK \langle d \rangle$ we put
$\bigl( U_t \bigr) _{ij} := j_{0t} (x_{ij} )$ to obtain a unitary operator $U_t$ on $\mc ^d \ot \SH$. Here $\SH$ is the completion of the pre-Hilbert space obtained through the Gelfand-Naimark-Segal construction applied to the pair $(\SA , \Phi)$. Indeed, the unitarity relations (\ref{unitaryRel}) force the operators $j_{st} (x_{ij})$ to be bounded and $U_t$ is a unitary matrix of bounded operators on $\SH$.
Moreover, $j _{st} (x_{ij} ) = \bigl( U_s ^{-1} U_t \bigr) _{ij}$, so that the process $j_{st}$ is determined by the \emph{unitary process} $U_t$. 
\par
One shows that all conditionally positive, hermitian linear functionals $\psi$ with $\psi (\ei ) = 0$  (that is all generators) on $\SK \langle d \rangle$ arise in the following way. (For a vector space $\SV$, denote by $\RM _d (\SV )$ the set of $d \times d$-matrices with entries in $\SV$.) Fix a pre-Hilbert space $H$ and a unitary operator $W$ on $\mc ^d \ot H$, a matrix $L$ in $\RM _d (H)$ and a hermitian matrix $D \in \RM _d (\mc )$. Then there is a unique generator $\psi$ such that $\rho (x_{ij} ) = W_{ij}$, 
$\eta (x_{ij} ) = L_{ij}$ and  $2 \ri \, D_{ij} = \psi (x_{ij} - x_{ji} ^* ) $.
The QSDE (\ref{QSDE}) for $U_t$ with generator triplet $(\rho , \eta , \psi )$ becomes
\be
\label{process}
\rd U_t = U_t \rd I _t , \  U_0 = \ei
\ee
with $I_t$ the matrix 
\begin{align}
\label{GNS}
\bigl(  I_t \bigr) _{ij} &= A_t \bigl( \eta (x_{ij} ^*) \bigr) + \GL _t \bigl( \rho (x_{ij}) - \gd _{ij} \bigr) + A _ t ^* \bigl( \eta (x_{ij} ) \bigr) + \psi (x_{ij} ) \, t  \\
&=
A_t ^* (L _ {ij} )  +  \GL _t  \bigl( (W - \ei )  _{ij} \bigr) - A_t  \bigl( (W ^* L)_{ji}  \bigr) + \bigl( D - \frac{1}{2}  L ^* L  \bigr) _{ij} \, t .
\end{align}
We have (cf. \cite{Fra06})
\begin{theorem}
All tensor L\'evy processes on $\SK \langle d \rangle$ are of type (\ref{process}).
\end{theorem}
For the process in \cite{vW84} we have $H = \mc$ and $W = \ei$, $D = 0$, so that (\ref{process}) can be written
\[
\rd U_t = U_t \bigl( L \rd A_t ^*  - L^* \rd A_t  - \frac{1}{2} L^*L  \rd t \bigr) .
\]

\section{Generalisations}

Again it was pioneering work of Wilhelm von Waldenfels that led to a new notion of independence, to \lq Boolean independence\rq . In the papers \cite{vW1} and \cite{vW75}, which again are motivated by physics, von Waldenfels used cumulants of Boolean independence at a time when there was no discussion at all on notions of non-commutative independence. Boolean independence had not even received its name. When later asked why he chose Boolean cumulants, he used to say: \lq Because they were the simplest to calculate\rq .
\bn
If we neglect questions of positivity, the category of classical probability spaces is formed by pairs $(\SA , \phi )$ as objects. Here $\SA$ is a commutative unital algebra and $\phi$ is a normalized linear functional on $\SA$, i.e. a linear map $\phi : \SA \to \mc$. The morphisms $j : (\SA _1 , \phi _1 ) \to (\SA _2 , \phi _2 )$ are given by algebra homomorphisms $j : \SA _1 \to \SA _2$ which satisfy $\phi _2 \circ j = \phi _1$ that is $\phi _1$ is the \lq distribution\rq \ of the \lq random variable\rq \ $j$. The joint distribution of two random variables $j_i : \SA _i  \to (\SA  , \phi )$, $i = 1, 2$, is the  linear functional $\phi \circ M _{\SA} \circ (j_1 \ot j _2 )$ on $\SA _1 \ot \SA _2 $. The random variables $j_1$, $j_2$ are called independent if their joint distribution equals $\phi _1 \ot \phi _2$, $\phi _i = \phi \circ j_i$ the distribution of $j_i$, $i = 1,2 $. This is the classical commutative situation. 
\par
One passes to non-commutativity by starting with the category $\SF$ with objects $(\SA , \phi )$ where $\SA$ is a (not necessarily unital) algebra and $\phi$ is a linear functional on $\SA$. (In order to include all known examples, for instance Boolean independence, we allow for non-unital algebras.)
Morphisms (= random variables) are defined as before as algebra homomorphisms with $\phi _1 = \phi _2 \circ j$. The joint distribution of $j_i$, $i = 1, 2$, is given as the linear functional $\phi \circ (j_1 \sqcup j_2 )$ on the \emph{free product} $\SA\sqcup \SA _2$ of $\SA _1 $ and $\SA _2$. The free product of algebras is the co-product in the category of algebras like the tensor product is  the co-product in the category of commutative unital algebras. $\SA _1 \sqcup \SA _2$ can be realized as the quotient of the tensor algebra $\RT (\SA _1 \oplus \SA _2 )$ by the ideal generated by the elements $a _ 1 \ot b_1 - a _1 b_1$, $a _1 , b_1 \in \SA _1$ and $a _ 2 \ot b_2 - a _2 b_2$, $a _2 , b_2 \in \SA _2$. 
\par
A \emph{universal product} is a bi-functor $\odot$ on  $\SF $ of the type $(\SA _1 , \phi _1 ) \odot (\SA _ 2 , \phi _2 ) = (\SA _1 \sqcup \SA _2 , \phi _1 \odot \phi _2 )$, $j_1 \odot j_2 = j_1 \sqcup j_2$,  such that $\odot$ turns $\SF$ into a \emph{monoidal category} with unit object $0  : \{ 0 \} \to \mc$; see 
\cite{Lachs1} and \cite{Malte}. The term \emph{universal} in this connection might sound strange for a category theorist. It comes from the condition that  the product \lq is the same\rq \ for all algebras which, of course, means nothing else but that it is the tensor product of a category. Given $\odot$ two random variables are called 
$\odot$-independent if their joint distribution equals $\phi _ 1\odot \phi _2$ that is the $\odot$-product of the marginal distributions. Examples of $\odot$ are  Boolean independence, tensor independence, free independence, monotone independence and
 anti-monotone independence. 
 \par
All these products are positive in the following sense. If $\SA _i$, $i = 1,2$, are $*$-algebras and if $\phi _i : \SA _i \to \mc$ are such that the \emph{unitizations} $\ei \phi _i$ are states on the unitizations $\ei \SA _i$, then the unitization $\ei
\bigl( \phi _1 \odot \phi _2 \bigr)$ is a state on $\ei \bigl( \SA _1 \sqcup \SA _2 \bigr)$, see \cite{Malte}. The first classification result for products $\odot$ was in the paper \lq On universal products\rq \ \cite{Spe97} by Roland Speicher. The complete classification of positive universal products was achieved by Naofumi Muraki in \cite{Mur03} where it is shown that the above mentioned five examples are the only products which satisfy an additional property (which is satisfied by positive products). For a classification of general universal products see \cite{GeLa}.
\bn
Still there are more examples. One arose by work of Marek Bo\. zejko and Roland Speicher who investigated \emph{conditional freeness}; see \cite{BSp}. More examples of this kind were found by Takahiro Hasebe; see \cite{Has}. In these examples the linear functional $\phi$ is generalized to a vector valued linear functional $\phi : \SA \to \mc ^d $, $d$ a fixed natural number. This gives many more possibilities and a full classification seems to be out of reach. 
\par
Another direction of generalization is to work with $m$-\emph{faced} algebras, $m \in \mn$,  i.e. with algebras $\SA$ which come with a free decomposition $\SA = \SA _1 \sqcup \ldots \sqcup \SA _m$. This means that the $\SA _i$ are sub-algebras of $\SA$ such that the natural mapping from $\SA _1 \sqcup \ldots \sqcup \SA _m $ to $\SA$ is an isomorphism of algebras. The $d$- and $m$-generalizations can be done in one step. Consider the category $\SF _{d, m}$ with objects pairs $( \SA , \phi )$ where $\SA$ is an $m$-faced algebra and $\phi : \SA \to \mc ^d $ linear. The morphisms $j$ are algebra homomorphisms which respect the free decompositions and with $\phi _1  = \phi _2 \circ j$. A natural product $\odot$, as before, is a  bi-functor on $\SF _{d, m}$ with  $(\SA _1 , \phi _1 ) \odot (\SA _ 2 , \phi _2 ) = (\SA _1 \sqcup \SA _2 , \phi _1 \odot \phi _2 )$, $j_1 \odot j_2 = j_1 \sqcup j_2$, and  such that $\odot$ turns $\SF$ into a monoidal category with unit object $0  : \{ 0 \} \to \mc ^d$; see again \cite{Malte}. $\odot$-independence is defined as before.
The first example of this kind was the bi-free independence which Dan V. Voiculescu introduced in 2014; see \cite{Voi123}. 
\par
Given such a product $\odot$ we are in the situation to define $\odot$-L\'evy processes on $m$-faced dual semi-groups. The latter are co-monoids in the monoidal category of $m$-faced algebras with product $\sqcup$ and unit object $\{ 0 \}$. More precisely, 
an $m$-faced dual semi-group is a triplet $(\SB , \GL , 0 )$ where $\SB$ is an $m$-faced $*$-algebra and $\GL : \SB \to \SB \sqcup \SB $ is an $m$-faced $*$-algebra homomorphism such that $(\GL \sqcup \id ) \circ \GL = (\id \sqcup \GL) \circ \GL$ and $(\id \sqcup 0 ) \circ \GL = \id = (0 \sqcup \GL ) \circ \GL$;  see \cite{Malte}. A $\odot$-L\'evy process on the $m$-faced dual semi-group $\SB$ over the object $(\SA , \phi )$ of $\SF _{d, m}$, is, in analogy to the tensor case, a family $j_{st} : \SB \to \SA$ of $m$-faced $*$-algebra homomorphisms such that the increment property, the stationarity and independence of increments, and the week continuity (which now reads $\Phi \circ j_{st}  \underset{t \downarrow 0}{\longrightarrow} 0$) are satisfied.
\par 
As before the 1-dimensional distributions $\phi _t = \Phi \circ j_{0, t }$ form a convolution semi-group where now the convolution is defined by 
\[
\phi _1 \star \phi _2 = (\phi _1 \odot \phi _2 ) \circ \GL .
\]
We also require our product $\odot$ to be positive. Using the \emph{Lachs functor}, one can show (see \cite{Lachs1} and \cite{Malte}) that $\frac{1}{t} \bigl( \phi _t (b)\bigr)$ converges to a limit which we denote by $\psi (b)$. Moreover, with the help of the Lachs functor $\phi _t$ can be reconstructed from $\psi$ as an \lq exponential\rq \  $\phi _t = \exp_{\star} \bigl( t \psi\bigr)$. A very nice result of Malte Gerhold (Theorem 5.9 in \cite{Malte}) says that in this situation the Schoenberg correspondence again holds. The proof of this general Schoenberg correspondence uses the tensor result (Theorem \ref{Schoenalg}) which had been proved in the co-commutative case by Wilhelm von Waldenfels. With this result we have that $\odot$-L\'evy processes are given by conditionally positive, hermitian linear functionals. 
\bn
We close by mentioning that a realization of $\odot$-L\'evy processes as solutions of QSDE is a nice problem. If $d = m = 1$ this, at least partially, has been done (free Fock space, Boolean Fock space). The general $d,m$-case seems to be an open problem.
A classification of products $\odot$ in the case $d = 1$, $m = 2$ is in progress and there exist very interesting partial results that can be found in the PhD thesis of Philipp Var\v{s}o \cite{Philipp}.

\bibliographystyle{alpha}
\bibliography{mybib-10,mybib-11}

\newcommand{\Swap}[2]{#2#1}\newcommand{\Sort}[1]{}
\begin{thebibliography}{ASvW88}

\bibitem[ASvW88]{waldenfels88}
L.~Accardi, M.~Sch\"urmann, and W.~von Waldenfels.
\newblock Quantum independent increment processes on superalgebras.
\newblock {\em Math. Z.}, 198:451--477, 1988.

\bibitem[Bro81]{Brown}
L.G. Brown.
\newblock Ext of certain free product {$C ^*$}-algebras.
\newblock {\em J. Operator Theory}, 6:135--141, 1981.

\bibitem[BS91]{BSp}
M.~Bo{\.z}ejko and R.~Speicher.
\newblock {$\psi$}-independent and symmetrized white noises.
\newblock In {\em Quantum probability \& related topics}, QP-PQ, VI, pages
  219--236. World Sci. Publ., River Edge, NJ, 1991.

\bibitem[CH71]{CuHu}
C.D. Cushen and R.L. Hudson.
\newblock A quantum-mechanical central limit theorem.
\newblock {\em Journal of Applied Probability}, 82:454--469, 1971.

\bibitem[CH77]{Hud}
A.M. Cockroft and R.L. Hudson.
\newblock Quantum mechanical {W}iener processes.
\newblock {\em J. Multivariate Anal.}, 7(1):107--124, 1977.

\bibitem[Chu68]{Chung}
K.L. Chung.
\newblock {\em A course in probability theory}.
\newblock Harcourt, Brace and World, 1968.

\bibitem[Fra06]{Fra06}
U.~Franz.
\newblock L\'evy processes on quantum groups and dual groups.
\newblock In {\em Quantum independent increment processes {II}}, volume 1866 of
  {\em Lecture Notes in Math.}, pages 161--257. Springer, Berlin, 2006.

\bibitem[Ger15]{Ger1}
M.~Gerhold.
\newblock {\em On several problems in the theory of comonoi\-dal sys\-tems and
  sub\-product sys\-tems}.
\newblock PhD thesis, Greifs\-wald, 2015.
\newblock \url{http://ub-ed.ub.uni-greifswald.de/opus/volltexte/2015/2244/}.

\bibitem[Ger21]{Malte}
M.~Gerhold.
\newblock Schoenberg correspondence for multifaced independences.
\newblock {\em arXiv:2104.02985v2}, 2021.

\bibitem[GL15]{GeLa}
M.~Gerhold and S.~Lachs.
\newblock Classification and {GNS}-construction for general universal products.
\newblock {\em Infin. Dimens. Anal. Quantum Probab. Relat. Top.},
  18(1):1550004, 29, 2015.

\bibitem[GvW78]{GivW}
N.~Giri and W.~von Waldenfels.
\newblock An algebraic central limit theorem in the anti-commuting case.
\newblock {\em Zeitschrift f{\"u}r Wahrscheinlichkeitstheorie und Verwandte
  Gebiete}, 42:129--134, 1978.

\bibitem[GvW89]{GlvW}
P.~Glockner and W.~von Waldenfels.
\newblock The relations of the non-commutative coefficient algebra of the
  unitary group.
\newblock In {\em Quantum Probability and Applications IV. Lecture Notes in
  Mathematics, vol 1396}. Springer, Berlin, Heidelberg, 1989.

\bibitem[Has11]{Has}
T.~Hasebe.
\newblock Conditionally monotone independence {I}: {I}ndependence, additive
  convolutions and related convolutions.
\newblock {\em Infin. Dimens. Anal. Quantum Probab. Relat. Top.},
  14(3):465--516, 2011.

\bibitem[HP84]{HudPar}
R.L. Hudson and K.R. Parthasarathy.
\newblock Quantum {I}to's formula and stochastic evolutions.
\newblock {\em Comm. Math. Phys.}, 93(3):301--323, 1984.

\bibitem[Lac15]{Lachs1}
S.~Lachs.
\newblock {\em A new family of universal products and aspects of a non-positive
  quantum probability theory}.
\newblock PhD thesis, Greifs\-wald, 2015.
\newblock \url{http://ub-ed.ub.uni-greifswald.de/opus/volltexte/2015/2242/}.

\bibitem[Mur03]{Mur03}
N.~Muraki.
\newblock The five independences as natural products.
\newblock {\em Infin. Dimens. Anal. Quantum Probab. Relat. Top.},
  6(3):337--371, 2003.

\bibitem[Par92]{Par92}
K.R. Parthasarathy.
\newblock {\em {An introduction to quan\-tum stochastic calculus}}.
\newblock Birk\-h\"au\-ser, 1992.

\bibitem[Sch91]{MSchue91b}
M.~Sch\"urmann.
\newblock {White noise on involutive bialgebras}.
\newblock In L.~Accardi, editor, {\em Quantum Probability \& Related Topics
  VI}, number~VI in QP-PQ:, pages 401--419. World Scientific, 1991.

\bibitem[Spe97]{Spe97}
R.~Speicher.
\newblock On universal products.
\newblock In {\em Free probability theory ({W}aterloo, {ON}, 1995)}, volume~12
  of {\em Fields Inst. Commun.}, pages 257--266. Amer. Math. Soc., Providence,
  RI, 1997.

\bibitem[Var]{Philipp}
Ph. Var{\v{s}}o.
\newblock {\em Studies on positive and symmetric two-faced universal products}.
\newblock PhD thesis, Greifs\-wald, thesis submitted 2021.

\bibitem[Voi14]{Voi123}
D.-V. Voiculescu.
\newblock Free probability for pairs of faces {I}.
\newblock {\em Comm. Math. Phys.}, 332(3):955--980, 2014.

\bibitem[vW73]{vW1}
W.~von Waldenfels.
\newblock {An approach to the theory of pressure broadening of spectral lines}.
\newblock In M.~Behara, K.~Krickeberg, and J.~Wolfowitz, editors, {\em
  Probability and Information Theory II}, number 296 in Lect.\ Notes Math.
  Springer, 1973.

\bibitem[vW75]{vW75}
W.~von Waldenfels.
\newblock Interval partitions and pair interactions.
\newblock In {\em Séminaire de probabilités (Strasbourg) Tome 9}, pages
  565--588. Springer, Berlin, 1975.

\bibitem[vW84a]{vW84}
W.~von Waldenfels.
\newblock It\^o solution of the linear quantum stochastic differential equation
  describing light emission and absorption.
\newblock In {\em Quantum probability and applications to the quantum theory of
  irreversible processes ({V}illa {M}ondragone, 1982)}, volume 1055 of {\em
  Lecture Notes in Math.}, pages 384--411. Springer, Berlin, 1984.

\bibitem[vW84b]{vW3}
W.~von Waldenfels.
\newblock Positive and conditionally positive sesquilinear forms on
  anticommutative coalgebras.
\newblock In {\em Probability measures on groups, {VII} ({O}berwolfach, 1983)},
  volume 1064 of {\em Lecture Notes in Math.}, pages 450--466. Springer,
  Berlin, 1984.

\end{thebibliography}

\end{document}